\documentclass{commat}

\title{On square-free numbers generated from given sets of primes}
\thanks{This is an updated version of the originally published article.}
\author{G\'abor Rom\'an}
\affiliation{\address{G\'abor Rom\'an -- E\"otv\"os Lor\'and University, Budapest, Hungary}\email{rogpaai@inf.elte.hu}}

\abstract{%
Let $x$ be a positive real number, and $\mathcal{P} \subset [2,\lambda(x)]$ be a set of primes,
where $\lambda(x) \in \Omega(x^\varepsilon)$ is a monotone increasing function with $\varepsilon \in (0,1)$.
We examine $Q_{\mathcal{P}}(x)$, where $Q_{\mathcal{P}}(x)$ is the element count of the set containing those positive square-free integers,
which are smaller than-, or equal to $x$, and which are only divisible by the elements of $\mathcal{P}$.
}



\VOLUME{30}
\NUMBER{1}
\firstpage{229}
\DOI{https://doi.org/10.46298/cm.10527}

\paper

\section{Introduction}
In this article, we are going to investigate the number of square-free numbers that one can generate from given sets of prime numbers.
More precisely, let's take a set of prime numbers $\mathcal{P}$,
and denote with $Q_{\mathcal{P}}(x)$ the element count of the set of all those positive square-free integers,
which are smaller than-, or equal to $x$; and which are only divisible by the elements of $\mathcal{P}$.
We would like to know how $Q_{\mathcal{P}}(x)$ behaves asymptotically based on the structure of $\mathcal{P}$.

Following the classical notation, denote the number of \textit{quadratfrei} integers between $1$ and $x$ as $Q(x)$.
We know that $Q(x) = x/\zeta(2) + \mathcal{O}(\sqrt{x})$ for all $x \geq 1$, see for example \cite[Th. 2.2]{MontgomeryVaughan2006}.
What is required from $\mathcal{P}$ to achieve $Q_{\mathcal{P}}(x) \asymp x$?
Indulge ourselves for a moment with an informal train of thoughts.
Taking only small primes, we can select even all of them to form products which aren't greater than $x$.
Then by the binomial theorem we have that the number of these products grow exponentially.
Based on this,
we can expect that there is a threshold for the number of primes in $\mathcal{P}$ somewhere between the (poly)logarithmic order,
and the fractional power order, where the behaviour of $Q_{\mathcal{P}}(x)$ changes.

%
%
\begin{proposition}
\label{proposition:count using all primes}
Let $\lambda : \mathbb{R} \rightarrow [1,+\infty)$ be a monotone increasing function which is in
$\Omega(x^\varepsilon)$ with $\varepsilon \in (0,1)$,
and let $\mathcal{P}$ contain all the primes which are not greater than $\lambda(x)$.
Then we have
\[
	Q_{\mathcal{P}}(x) \asymp x\frac{\ln\lambda(x)}{\ln{x}}
\]
as $x \rightarrow \infty$.
\end{proposition}

\section{Proofs}
We are going to use the fundamental lemma of the combinatorial sieve to prove our results, see \cite[P I.4, Th. 4.4]{Tenenbaum2015}.
The idea is to utilise the sieve to remove those square-free integers below $x$, which are not divisible by a prime in $\mathcal{P}$,
so we will sift with the set $\mathcal{P}^\prime := [1,x]\cap\mathbb{P} \setminus \mathcal{P}$.
Using the notation of \cite[P I.4, Th. 4.4]{Tenenbaum2015}, let $\mathcal{A}$ be the set of all square-free numbers not greater than $x$,
let $\mathcal{A}_d := \lbrace{a \in \mathcal{A} : a \equiv 0\ (\textrm{mod}\ d)}\rbrace$, and finally let
\[
	P(y) := \prod_{\substack{p \leq y \\ p \in \mathcal{P}^\prime}} p.
\]
We are going to give an asymptotic for
\[
	S(\mathcal{A}, \mathcal{P}^\prime, x) := |\lbrace{a \in \mathcal{A} : (a,P(x)) = 1}\rbrace|
\]
using the mentioned sieve.
To do this, on the one hand we need a multiplicative function $w \geq 0$ such that
\[
	|\mathcal{A}_d| = \frac{w(d)}{d}X + R_d
\]
for some real number $X$, and all $d|P(x)$.
We are going to give an asymptotic for $|\mathcal{A}_d|$ in section \ref{subsection:Cardinality of the sieving sets},
and based on it an appropriate multiplicative function.
On the other hand, we need to show that there exist positive constants $\kappa$ and $\kappa^\prime$ such that the inequality
\begin{align}
\label{align:requirement 2}
	\prod_{\eta \leq p \leq \xi} \bigg(1 - \frac{w(p)}{p}\bigg)^{-1}
	< \bigg(\frac{\ln\xi}{\ln\eta}\bigg)^{\kappa} \bigg(1 + \frac{\kappa^\prime}{\ln\eta}\bigg)
\end{align}
holds for $2 \leq \eta \leq \xi$ using our multiplicative function $w$.
We are going to look at this in section \ref{subsection:The Euler product}.
If these requirements hold, then uniformly in $\mathcal{A}$, $X$, and $u \geq 1$ we have
\begin{align}
\label{align:fundamental lemma of the combinatorial sieve}
	S(\mathcal{A}, \mathcal{P}^\prime, x)
	= X \prod_{\substack{p \leq x \\ p \in \mathcal{P}^\prime}} \bigg(1 - \frac{w(p)}{p}\bigg) (1 + \mathcal{O}(u^{-u/2}))
	+ \mathcal{O}\bigg(\sum_{\substack{d \leq x^u \\ d|P(x)}} |R_d|\bigg)
\end{align}
as $x \rightarrow +\infty$, using \cite[P I.4, Th. 4.4]{Tenenbaum2015}.
We are going to use this equation to prove our results in section \ref{subsection:The main asymptotic}.

We will rely on the Dedekind psi function
\begin{align}
\label{align:Dedekind psi function}
	\psi(d) = d \prod_{p|d} \bigg(1 + \frac{1}{p}\bigg)
\end{align}
see \cite[B41]{Guy2004}, or \cite[Ch. 3, Ex. 11]{Apostol1976};
and the following representation of the Euler totient function
\begin{align}
\label{align:Euler totient function}
	\varphi(d) = d \prod_{p|d} \bigg(1 - \frac{1}{p}\bigg)
\end{align}
see \cite[Th. 2.4]{Apostol1976}.
Recall that the equality
\begin{align}
\label{align:Euler totient function, sum representation}
	\frac{\varphi(d)}{d} = \sum_{q|d} \frac{\mu(q)}{q}
\end{align}
holds when $d \geq 1$, see \cite[Th. 2.3]{Apostol1976}.

\subsection{Cardinality of the sieving sets}
\label{subsection:Cardinality of the sieving sets}

%
%
\begin{lemma}
\label{lemma:cardinality of the sieving sets}
For every $\delta > 0$, and $d \in Q(x^{2/3-\eta})$ with $\eta > 0$ we have
\begin{align}
\label{align:cardinality of the sieving sets}
	|\mathcal{A}_d| = \frac{x}{\zeta(2)\psi(d)} + \mathcal{O}(x^{1/2}d^{-1/4}) + \mathcal{O}(d^{1/2+\delta})
\end{align}
as $x \rightarrow +\infty$.
\end{lemma}
To prove this lemma we are going to follow the method of Prachar, see article \cite{Prachar1958}.
Take note that the error term could be improved, see the article of Hooley \cite{Hooley1975}.
For further developments see the article of Nunes \cite{Nunes2015} and the article of Mangerel \cite{Mangerel2021}.

%
%
\begin{proof}[Proof of lemma \ref{lemma:cardinality of the sieving sets}]
Fix a big enough real $x > 0$, and choose an appropriate $d$. Then
\[
	|\mathcal{A}_d|
	= \sum_{\substack{n \leq x \\ n \equiv 0 (d)}} \mu(n)^2
	= \sum_{\substack{n \leq x \\ n \equiv 0 (d)}} \sum_{l^2|n} \mu(l)
	= \sum_{\substack{l^2m \leq x \\ l^2m \equiv 0 (d)}} \mu(l)
	= \sum_{l \leq X_d} + \sum_{X_d < l}
\]
where $X_d := x^{1/2}d^{-1/4}$. We look at the two sums separately.
\begin{itemize}
\item[$\sum_{l \leq X_d}$]
The congruence $l^2m \equiv 0\ (\mathrm{mod}\ d)$ is always soluble for $m$ because $(l^2,d)|0$;
and the solutions can be given as $m \equiv 0\ (\mathrm{mod}\ d/(l,d))$ because we have $(l^2,d) = (l,d)$ as $d$ is square-free.
Thus the first sum is
\[
	\sum_{l \leq X_d} \mu(l) \sum_{\substack{m \leq x/l^2 \\ m \equiv 0 (d/(l,d))}} 1 
	= \sum_{l \leq X_d} \mu(l) \bigg(\frac{(l,d)x}{l^2d} + \mathcal{O}(1)\bigg)
\]
where the second item between the parentheses contributes $\mathcal{O}(X_d)$ to the result, so we are going to focus on the first item.
We can split the sum along the divisors of $d$ as
\[
	\frac{x}{d} \sum_{l \leq X_d} \frac{(l,d)\mu(l)}{l^2}
	= \frac{x}{d} \sum_{q|d} \sum_{\substack{l \leq X_d \\ (l,d) = q}} \frac{q\mu(l)}{l^2}.
\]
For a given $q$, we can write the appropriate $l$ as $qm$ with $(q,m)=1$.
(Observe that we only have to worry about the square-free $l$.)
In this case $q=(qm,d)=(q,d)(m,d)$, so $(m,d)=1$, and
\[
	\frac{x}{d} \sum_{q|d} \sum_{\substack{m \leq X_d/q \\ (m,d) = 1}} \frac{q\mu(qm)}{(qm)^2}
	= \frac{x}{d} \sum_{q|d} \frac{\mu(q)}{q} \sum_{\substack{m \leq X_d/q \\ (m,d) = 1}} \frac{\mu(m)}{m^2}
\]
as $(q,m)=1$, and $\mu$ is multiplicative.
The internal sum on the right hand side can be written as
\[
	\sum_{\substack{m=1 \\ (m,d)=1}}^\infty \frac{\mu(m)}{m^2} + \mathcal{O}\bigg(\sum_{m > X_d/q} \frac{1}{m^2}\bigg)
	= \frac{1}{\zeta(2)} \prod_{p|d}\bigg(1 - \frac{1}{p^2}\bigg)^{-1} + \mathcal{O}(q/X_d).
\]
Relying on (\ref{align:Euler totient function, sum representation}) we get
\begin{align}
	\frac{1}{d} \sum_{q|d} \frac{\mu(q)}{q} \prod_{p|d}\bigg(1 - \frac{1}{p^2}\bigg)^{-1}
	&= \frac{\varphi(d)}{d^2} \prod_{p|d}\bigg(1 - \frac{1}{p^2}\bigg)^{-1} \nonumber \\
	&= \frac{1}{d} \prod_{p|d}\bigg(1 + \frac{1}{p}\bigg)^{-1}
	= \frac{1}{\psi(d)} \nonumber
\end{align}
by using (\ref{align:Dedekind psi function}), and (\ref{align:Euler totient function}).
Concerning the sum of the remaining asymptotic term, we have
\[
	\frac{x}{d} \sum_{q|d} \frac{\mu(q)}{q} \mathcal{O}(q/X_d)
	= \frac{x}{dX_d} \sum_{q|d} c_{d,q}\mu(q)
	\in \mathcal{O}(X_d)
\]
where $c_{d,q} \in \mathcal{O}(1)$.
So we have
\[
	\sum_{l \leq X_d} = \frac{x}{\zeta(2)\psi(d)} + \mathcal{O}(X_d).
\]
\item[$\sum_{X_d < l}$]
Observe that, since $X_d < l$, we have $m < d^{1/2}$; and thus, $l \leq (x/m)^{1/2}$ for a fixed $m$.
Taking a fixed $m$, if the congruence $l^2m \equiv 0\ (\mathrm{mod}\ d)$ is soluble, then the number of solutions is $2^{\omega(d)}$,
see \cite[Th. 5.28]{Apostol1976}, where $\omega(d)$ denotes the number of distinct prime factors of $d$,
with the convention that $\omega(1) = 0$.
Once the solutions are found, we can take all the admissible values of $l$, for which $l^2$ is congruent to a solution $\mathrm{mod}\ d$.
So we have that the second sum can be bounded by some constant times
\[
	\sum_{m < d^{1/2}} 2^{\omega(d)} \bigg(\bigg(\frac{x}{m}\bigg)^{1/2}d^{-1} + \mathcal{O}(1)\bigg)
	\ll 2^{\omega(d)}(x^{1/2}d^{-3/4} + d^{1/2})
\]
which is in $\mathcal{O}(x^{1/2}d^{-3/4+\delta}) + \mathcal{O}(d^{1/2+\delta})$ because $\omega(d) \in \mathcal{O}((\ln{d})/\ln_2{d})$, see \cite[Sec.~5.3]{Tenenbaum2015}.
\end{itemize}
Combining the two results, we get our statement.
\end{proof}
%
%

%
%
\begin{lemma}
\label{lemma:cardinality of the sieving sets, selberg-delange}
For $d \in Q(x)$ we have
\[
	|\mathcal{A}_d| = \frac{x}{\zeta(2)\psi(d)} + \mathcal{O}\bigg(\frac{x}{d}e^{-c\sqrt{\ln{x/d}}}\bigg)
\]
as $x \rightarrow +\infty$, where $c$ is a positive constant.
\end{lemma}

%
%
\begin{proof}[Proof of lemma \ref{lemma:cardinality of the sieving sets, selberg-delange}]
Let $d \in Q(x)$. We can write
\[
	|\mathcal{A}_d|
	= \sum_{\substack{n \leq x \\ n \equiv 0 (d)}} \mu(n)^2
	= \sum_{m \leq x/d} \mu(md)^2.
\]
The summand is zero when $(m,d) > 1$, otherwise we have $\mu(md)^2 = \mu(m)^2$ because $\mu$ is multiplicative, and $d$ is square-free.
So we have
\begin{align}
\label{align:cardinality of the sieving sets, selberg-delange 1}
	\sum_{\substack{m \leq x/d \\ (m,d)=1}} \mu(m)^2
	= \sum_{m \leq x/d} \chi_d(m) \mu(m)^2
\end{align}
where $\chi_d$ is the principal character modulo $d$.
Let
\[
	F_d(s;z)
	:= \sum_{n \geq 1} \frac{\chi_d(n) \mu(n)^2 z^{\omega(n)}}{n^s}.
\]
The function in the numerator of the summand is multiplicative, and when $s > 1$ we have
\[
	\sum_p \sum_{\nu \geq 1} \bigg|\frac{\chi_d(p^\nu) \mu(p^\nu)^2 z^{\omega(p^\nu)}}{p^{\nu s}}\bigg|
	= |z| \sum_p \bigg|\frac{\chi_d(p)}{p^s}\bigg|
	< +\infty
\]
so we can apply \cite[P. II.1, Th. 1.3]{Tenenbaum2015} to get that $F_d$ is absolutely convergent, and
\[
	F_d(s;z)
	= \prod_p \sum_{\nu \geq 0} \frac{\chi_d(p^\nu) \mu(p^\nu)^2 z^{\omega(p^\nu)}}{p^{\nu s}}
	= \prod_p \bigg(1 + \frac{\chi_d(p)z}{p^s}\bigg)
\]
when $s > 1$.
Now we are going to rely on the technique of \cite[P. II.6, Sec. 6.1]{Tenenbaum2015}.
The function
\[
	G_d(s; z)
	:= F_d(s; z) \zeta(s)^{-z}
	= \prod_p \bigg(1 + \frac{\chi_d(p)z}{p^s}\bigg) \bigg(1 - \frac{1}{p^s}\bigg)^z.
\]
is expandable as a Dirichlet series $G_d(s; z) = \sum_{n \geq 1} g_{d,z}(n)/n^s$,
where $g_{d,z}$ is the multiplicative function whose values on prime powers are determined by the following identities.
\begin{itemize}
\item
When $(p,d) = 1$, then $\chi_d(p) = 1$; and
\[
	1 + \sum_{\nu \geq 1} g_{d,z}(p^\nu) \xi^\nu
	= (1 + \xi z)(1 - \xi)^z
\]
where $|\xi| < 1$.
Using the binomial theorem on the right hand side we get
\begin{align}
	(1 + \xi z) \sum_{k=0}^\infty \binom{z}{k} (-\xi)^k
	&= (1 + \xi z)\bigg(1 - \xi z + \binom{z}{2} \xi^2 - \ldots\bigg) \nonumber \\
	&= 1 + \xi^2 \bigg(\binom{z}{2} - z^2\bigg) + \ldots \nonumber
\end{align}
so $g_{d,z}(p) = 0$.
Using Cauchy's inequality, see for example \cite[Sec. 2.5]{Titchmarsh1939}, we get for $|z| \leq A$ that
\[
	|g_{d,z}(p^\nu)|
	\leq M 2^{\nu/2}
\]
for $\nu \geq 2$, where
\[
	M
	:= \sup_{\substack{|z| \leq A \\ |\xi| \leq 1/\sqrt{2}}} |(1 + \xi z)(1 - \xi)^z|.
\]
Thus when $(p,d) = 1$ we have
\[
	\sum_{\nu \geq 1} \frac{|g_{d,z}(p^\nu)|}{p^{\nu \sigma}}
	\leq M \sum_{\nu \geq 2} \frac{2^{\nu/2}}{p^{\nu \sigma}}
	= 2M \frac{1}{p^\sigma(p^\sigma - \sqrt{2})}.
\]
\item
When $(p,d) > 1$, then $\chi_d(p) = 0$; and
\[
	1 + \sum_{\nu \geq 1} g_{d,z}(p^\nu) \xi^\nu
	= (1 - \xi)^z
	= 1 - \binom{z}{1}\xi + \binom{z}{2}\xi^2 - \binom{z}{3}\xi^3 + \ldots
\]
holds.
This means that
\[
	\sum_{\nu \geq 1} \frac{|g_{d,z}(p^\nu)|}{p^{\nu \sigma}}
	= \sum_{\nu \geq 1} \binom{z}{\nu} \frac{1}{p^{\nu\sigma}}
	= \bigg(1 + \frac{1}{p^\sigma}\bigg)^z - 1.
\]
\end{itemize}

Now that we've given how $g_{d,z}$ behaves on prime powers, let's look at the absolute convergence of $G_d(s;z)$.
We have
\[
	\sum_{p} \sum_{\nu \geq 1} \frac{|g_{d,z}(p^\nu)|}{p^{\nu \sigma}}
	= \sum_{(p,d)=1} \sum_{\nu \geq 1} \frac{|g_{d,z}(p^\nu)|}{p^{\nu \sigma}}
	+ \sum_{(p,d)>1} \sum_{\nu \geq 1} \frac{|g_{d,z}(p^\nu)|}{p^{\nu \sigma}}
\]
where we look at the two sums separately.
\begin{itemize}
\item
For $\sigma > 1/2$ we have that the first sum is less than or equal to
\[
	2M \sum_{(p,d)=1} \frac{1}{p^\sigma(p^\sigma - \sqrt{2})}
	\leq 2M \sum_p \frac{1}{p^\sigma(p^\sigma - \sqrt{2})}
	\leq \frac{cM}{\sigma - 1/2}
\]
where $c$ is an absolute constant.

\item
The second sum is equal to
\[
	\sum_{(p,d)>1} \bigg\lbrace\bigg(1 + \frac{1}{p^\sigma}\bigg)^z - 1\bigg\rbrace
\]
which is finite as $d$ has finite number of prime divisors.
\end{itemize}
It follows that $G_d(s;z)$ is absolutely convergent for $\sigma > 1/2$, and for $\sigma \geq 3/4$, we have $G_d(s;z) \ll_{A,d} 1$.
Based on \cite[P. II.5, Th. 5.2]{Tenenbaum2015}, for $x \geq 3$, $N \geq 0$, $A > 0$, $0 < z \leq A$, we have
\[
	\sum_{n \leq x} \chi_d(n) \mu(n)^2 z^{\omega(n)}
	= x (\ln{x})^{z-1} \bigg\lbrace{\sum_{0 \leq k \leq N} \frac{\lambda_k(z)}{(\ln{x})^k} + \mathcal{O}(R_N(x))}\bigg\rbrace
\]
where
\[
	\lambda_k(z)
	:= \frac{1}{\Gamma(z-k)} \sum_{h+j=k} \frac{1}{h!j!} G_d^{(h)}(1;z) \gamma_j(z)
\]
with $\gamma_j(z)$ defined as in \cite[P. II.5, Th. 5.1]{Tenenbaum2015}, and $R_N(x) := R_N(x, c_1, c_2)$.
The positive constants $c_1$, $c_2$, and the implicit constant in the Landau symbol depend at most on $A$ and the convergence properties of $G_d$.

We are interested in the case when $z=1$.
As $\lambda_k(z) = 0$ whenever $k \geq z$, we only have to deal with $\lambda_0(1)$.
We can also choose $N$ to minimise the error term, see the discussion after \cite[P. II.5, Th. 5.2]{Tenenbaum2015},
so at the end of the day we have
\begin{align}
\label{align:cardinality of the sieving sets, selberg-delange 2}
	\sum_{n \leq x} \chi_d(n) \mu(n)^2
	= x \bigg\lbrace{\lambda_0(1) + \mathcal{O}(e^{-c\sqrt{\ln{x}}})}\bigg\rbrace
\end{align}
where $c$ is a positive constant.
As $\gamma_0(1) = 1$ holds, see \cite[P. II.5, Th. 5.1]{Tenenbaum2015},
we have $\lambda_0(1) = G_d(1;1)$ which is
\[
	\prod_p \bigg(1 + \frac{\chi_d(p)}{p}\bigg) \bigg(1 - \frac{1}{p}\bigg)
	= \prod_p \bigg(1 - \frac{1}{p^2}\bigg) \prod_p \bigg(1 + \frac{1 - \chi_d(p)}{p}\bigg)^{-1}
	= \frac{1}{\zeta(2)}	\prod_{p|d} \bigg(1 + \frac{1}{p}\bigg)^{-1}
\]
as
\[
	\prod_p \bigg(1 + \frac{\chi_d(p)}{p}\bigg)
	= \prod_p \bigg(1 + \frac{1}{p}\bigg) \prod_p \bigg(1 + \frac{1 - \chi_d(p)}{p}\bigg)^{-1}.
\]
We get the desired result from (\ref{align:cardinality of the sieving sets, selberg-delange 1})
using (\ref{align:cardinality of the sieving sets, selberg-delange 2}).
\end{proof}
%
%

We can set our multiplicative function as
\begin{align}
\label{align:multiplicative function}
	w(d) := \prod_{p|d}\bigg(1+\frac{1}{p}\bigg)^{-1}
\end{align}
and $X$ as $x/\zeta(2)$, furthermore choose $R_d$ based on $d$ to have the required form for $|\mathcal{A}_d|$.

\subsection{The Euler product}
\label{subsection:The Euler product}

%
%
\begin{lemma}
\label{lemma:euler product}
For every $y \geq 2$, we have that the inequalities
\[
	\frac{4}{3} \frac{e^{-\gamma}}{\ln{y}}\bigg(1 - \frac{1}{\ln^2{y}}\bigg)
	< \prod_{p \leq y} \bigg(1-\frac{1}{p+1}\bigg) <
	\frac{\pi^2}{6} \frac{e^{-\gamma}}{\ln{y}}\bigg(1 + \frac{1}{\ln^2{y}}\bigg)
\]
hold, where $\gamma$ is the Euler-Mascheroni constant.
\end{lemma}
Note that we could improve the exponent of the logarithmic part inside the parentheses on the right hand side,
but this form will suffice for our needs.

%
%
\begin{proof}
Observe that
\begin{align}
\label{align:euler product 1}
	\prod_{p \leq y} \bigg(1-\frac{1}{p+1}\bigg) \prod_{p \leq y} \bigg(1-\frac{1}{p}\bigg)^{-1}
	= \prod_{p \leq y} \bigg(1-\frac{1}{p^2}\bigg)^{-1}
\end{align}
where we have the partial Euler product of the Riemann zeta function on the right hand side, see \cite[Sec. 11.5]{Apostol1976}.
Take note that as $\zeta(2) = \pi^2/6$, we have that the value of this product will be in $[4/3,\pi^2/6)$.
Based on \cite[Th. 7, Col.]{RosserSchoenfeld1962}, we have
\begin{align}
\label{align:euler product 2}
	\frac{e^{-\gamma}}{\ln{y}}\bigg(1 - \frac{1}{\ln^2{y}}\bigg)
	< \prod_{p \leq y} \bigg(1-\frac{1}{p}\bigg) <
	\frac{e^{-\gamma}}{\ln{y}}\bigg(1 + \frac{1}{\ln^2{y}}\bigg)
\end{align}
for every $y > 1$.
Combining (\ref{align:euler product 1}) and (\ref{align:euler product 2}) we get our statement.
\end{proof}

Now we are going to show that there exist positive constants $\kappa$,
and $\kappa^\prime$ for our function $w$ so that requirement (\ref{align:requirement 2}) holds.
Take note that the left hand side of requirement (\ref{align:requirement 2}) is
\begin{align}
\label{align:euler product}
	\prod_{\eta \leq p \leq \xi} \bigg(1 - \frac{1}{p+1}\bigg)^{-1}
\end{align}
when we substitute $w$, see (\ref{align:multiplicative function}).
\begin{itemize}
\item[$\eta = 2$]
When $2 \leq \xi < 3$, then the product will be $3/2$ so for example a $\kappa \geq 1$, and a $\kappa^\prime \geq 1$ will suffice.
Otherwise, when $\xi \geq 3$, then we can write
\[
	\prod_{p \leq \xi} \bigg(1 - \frac{1}{p+1}\bigg)^{-1}
	< (\ln\xi) \frac{3e^{\gamma}}{4} \bigg(1 - \frac{1}{\ln^2\xi}\bigg)^{-1}
	< 8 \ln\xi
\]
so a $\kappa \geq 1$, and a $\kappa^\prime \geq 5$ will suffice in this case.

\item[$2 < \eta < 3$]
When $\eta \leq \xi < 3$, then the product will be $1$, for which the previous constants will be good.
Otherwise, when $\xi \geq 3$, then
\[
	\prod_{2 < p \leq \xi} \bigg(1 - \frac{1}{p+1}\bigg)^{-1}
	< (\ln\xi) \frac{e^{\gamma}}{2} \bigg(1 - \frac{1}{\ln^2\xi}\bigg)^{-1}
	< 6 \ln\xi
\]
so a $\kappa \geq 1$, and a $\kappa^\prime \geq 6$ will suffice in this case.

\item[$\eta \geq 3$]
We can write product (\ref{align:euler product}) as
\[
	\prod_{p \leq \eta} \bigg(1 - \frac{1}{p+1}\bigg) \prod_{p \leq \xi} \bigg(1 - \frac{1}{p+1}\bigg)^{-1}
	< \bigg(\frac{\ln\xi}{\ln\eta}\bigg)^3 \frac{\pi^2}{8} \frac{(\ln^2\eta) + 1}{(\ln^2\xi) - 1}
\]
which shows us that a $\kappa \geq 3$ will suffice, and where
\[
	\frac{(\ln^2\eta) + 1}{(\ln^2\xi) - 1}
	< \frac{(\ln^2\eta) + 1}{(\ln^2\eta) - 1}
	= \frac{2}{(\ln^2\eta)-1} + 1
\]
thus we have to guarantee that the inequality
\[
	\frac{\pi^2}{8} \bigg(\frac{2}{(\ln^2\eta)-1} + 1\bigg) < 1 + \frac{\kappa^\prime}{\ln\eta}
\]
holds, which can be done by selecting a $\kappa^\prime \geq 14$.
\end{itemize}
Taking all the cases into consideration, it can be seen that we can select the required positive constants $\kappa$, and $\kappa^\prime$.

\subsection{The main asymptotic}
\label{subsection:The main asymptotic}
Select a function $\lambda$ satisfying the requirements of proposition \ref{proposition:count using all primes}.
$\mathcal{P}$ contains all the primes which are not greater than $\lambda(x)$,
so we have that the product which we have to compute in expression (\ref{align:fundamental lemma of the combinatorial sieve}) is equal to
\[
	\prod_{p \leq x} \bigg(1 - \frac{1}{p+1}\bigg) \prod_{p \leq \lambda(x)} \bigg(1 - \frac{1}{p+1}\bigg)^{-1}
	\asymp \frac{\ln\lambda(x)}{\ln{x}}
\]
based on lemma \ref{lemma:euler product}.

Because $|\mathcal{A}_d| = 0$ when $d > x$, we also have that $R_d = 0$ in this case.
So we can choose $u$ arbitrarily large,
rendering it ineffective in the product of expression (\ref{align:fundamental lemma of the combinatorial sieve}),
while in the error term it suffice to sum just until $x$.
By fixing a $\delta > 0$, we can write the error term as
\begin{align}
\label{align:error term}
	\sum_{\substack{d \leq x^{2/3 - \eta} \\ d|P(x)}} |R_d| + \sum_{\substack{x^{2/3 - \eta} < d \leq x \\ d|P(x)}} |R_d|
\end{align}
with an $\eta > \delta$, see lemma \ref{lemma:cardinality of the sieving sets}.
Then for big enough $x$, we can use equality (\ref{align:cardinality of the sieving sets}) to write the first sum as
\[
	c_1 x^{1/2} \sum_{\substack{d \leq x^{2/3 - \eta} \\ d|P(x)}} d^{-1/4} +
	c_2 \sum_{\substack{d \leq x^{2/3 - \eta} \\ d|P(x)}} d^{1/2+\delta}
\]
which is in $O(x^{1-\varepsilon(\delta,\eta)})$ with $\varepsilon(\delta,\eta) > 0$,
even if we don't take the structure of $\mathcal{P}$ into consideration.

What remains is to handle the second sum in expression (\ref{align:error term}).
We can use lemma \ref{lemma:cardinality of the sieving sets, selberg-delange} to write this sum as
\begin{align}
\label{align:error term split}
	c_3 x \sum_{\substack{x^{2/3 - \eta} < d \leq x \\ d|P(x)}} d^{-1} e^{-c\sqrt{\ln{x/d}}}
	= c_3 x \bigg\lbrace
	\sum_{\substack{x^{2/3 - \eta} < d \leq x^\varepsilon \\ d|P(x)}} d^{-1} e^{-c\sqrt{\ln{x/d}}}
	+ \sum_{\substack{x^\varepsilon < d \leq x \\ d|P(x)}} d^{-1} e^{-c\sqrt{\ln{x/d}}}
	\bigg\rbrace.
\end{align}
If $\varepsilon  < 2/3 - \eta$, then the first sum on the right hand side disappears.
Define 
\[
	A(y)
	:= \sum_{\substack{n \leq y \\ n|P(x)}} 1.
\]
If $y < x^\varepsilon$, then $A(y) = 0$; otherwise when $x^\varepsilon \leq y \leq x$, then
\[
	\pi(y) - \pi(x^\varepsilon)
	\leq A(y) \leq
	\Phi(y, x^\varepsilon)
\]
where $\Phi(y,z)$ is the count of those $n \leq y$ for which $P^-(n) > z$,
with $P^-(n)$ denoting the smallest prime factor of $n$, see \cite[P. III.6]{Tenenbaum2015}.
We have
\begin{align}
\label{align:phi}
	\Phi(y,z)
	\ll \frac{y}{\ln{z}}
\end{align}
for $2 \leq z \leq y$, see \cite[P. III.6, Th. 6.2]{Tenenbaum2015}.
Based on equation (\ref{align:phi}),
and on the known approximations for the prime counting function, see \cite{RosserSchoenfeld1962}, we can conclude that
\[
	A(y)
	\asymp \frac{y}{\ln{x}}
\]
when $2 \leq x^\varepsilon \leq y \leq x$.
Now we can use Abel summation to compute the sums on the right hand side of (\ref{align:error term split}).
For the first sum, we get
\[
	\frac{A(x^\varepsilon)}{x^\varepsilon}
	- \int_{x^{2/3 - \eta}}^{x^\varepsilon} A(t) \frac{d}{dt} \frac{e^{-c\sqrt{\ln{x/t}}}}{t} \,dt
	\asymp \frac{1}{\ln{x}}
\]
and for the second sum we get
\[
	\frac{A(x)}{x}
	- \frac{A(x^\varepsilon)}{x^\varepsilon} e^{-c\sqrt{\ln{x^{1 - \varepsilon}}}}
	- \int_{x^\varepsilon}^x A(t) \frac{d}{dt} \frac{e^{-c\sqrt{\ln{x/t}}}}{t} \,dt
	\asymp \frac{1}{\ln{x}}
\]
for big enough $x$ in both cases.

\section{Acknowledgments}

The author wishes to thank Tomohiro Yamada for pointing out an error in the argument of the article,
furthermore Attila B\'erczes and Ivan Kaygorodov for their support.


\EditInfo{July 29, 2021}{April 12, 2022}{Attila Bérczes}

\end{paper}